\documentclass[11pt,reqno, twoside]{amsart}
\usepackage{amscd}
\usepackage{amsfonts}
\usepackage{amsmath}
\usepackage{amssymb}
\usepackage{amsthm}
\usepackage{fancyhdr}
\usepackage{latexsym}
\usepackage[colorlinks=true, pdfstartview=FitV, linkcolor=blue, citecolor=blue, urlcolor=blue]{hyperref}
\usepackage{enumitem} 
\usepackage{mathrsfs}     
\usepackage{mathtools}
\usepackage{upgreek}
\usepackage{indentfirst} 
\usepackage{schemata} 
\usepackage{blindtext, rotating}   
\usepackage{soul} 
\setstcolor{red}
\usepackage{color}
\usepackage{tikz}
\usetikzlibrary{arrows.meta}
\usetikzlibrary{decorations.markings}
\tikzset{->-/.style={decoration={
  markings,
  mark=at position #1 with {\arrow{>}}},postaction={decorate}}}
  \tikzset{middlearrow/.style={
        decoration={markings,
            mark= at position 0.55 with {\arrow{#1}} ,
        },
        postaction={decorate}
    }
}
\usepackage{caption}           
\newcommand{\n}[2]{{\left\| #1 \right\|}_{#2}}
\newcommand{\f}[2]{\frac{#1}{#2}}
\newcommand{\lan}[1]{\left\langle #1 \right \rangle}
\newcommand{\ve}{\varepsilon}
\newcommand{\al}{\alpha}

\newcommand{\eee}[1]{\begin{equation}#1\end{equation}}
\newcommand{\sss}[1]{\begin{subequations}#1\end{subequations}}

\newcommand{\ddd}[1]{\begin{alignat}{2}#1\end{alignat}}

\newcommand{\nn}{\nonumber}
\newcommand{\p}{\partial}

\definecolor{ddgreen}{RGB}{0,170,0}

\renewcommand{\b}{\mathcolor{blue}}

\newcommand{\no}[1]{\left\| #1 \right\|}

\def\XXint#1#2#3{{\setbox0=\hbox{$#1{#2#3}{\int}$ }
\vcenter{\hbox{$#2#3$ }}\kern-.6\wd0}}

\makeatletter
\renewcommand\subsection{\@startsection{subsection}{2}%
  \z@{-.5\linespacing\@plus-.7\linespacing}{.5\linespacing}%
  {\normalfont\bfseries}}
\renewcommand\subsubsection{\@startsection{subsubsection}{3}%
  \z@{.5\linespacing\@plus.7\linespacing}{-.5em}%
  {\normalfont\scshape}}
\makeatother
\makeatletter
\def\mathcolor#1#{\@mathcolor{#1}}
\def\@mathcolor#1#2#3{%
\protect\leavevmode
\begingroup
\color#1{#2}#3%
\endgroup
}
\makeatother
\theoremstyle{plain}  
\newtheorem{theorem}{Theorem}[section]
\newtheorem{proposition}{Proposition}[section]
\newtheorem{lemma}{Lemma}[section]

\theoremstyle{definition}

\newtheorem{remark}{Remark}[section]
\newenvironment{Proof}[1][\proofname]
{\proof[\textnormal{\textbf{#1.}}]}{\endproof}
\newcommand{\bp}{\begin{Proof}}
\newcommand{\ep}{\end{Proof}}

\DeclareMathSizes{12}{12}{7}{5}
\numberwithin{figure}{section}
\numberwithin{equation}{section}
\usepackage{geometry}
\makeatletter
\def\l@section{\@tocline{1}{0pt}{1pc}{}{}}
\def\l@subsection{\@tocline{2}{0pt}{1pc}{4.6em}{}}
\def\l@subsubsection{\@tocline{3}{0pt}{1pc}{7.6em}{}}
\renewcommand{\tocsection}[3]{%
  \indentlabel{\@ifnotempty{#2}{\makebox[2.3em][l]{%
    \ignorespaces#1 #2.\hfill}}}#3}
\renewcommand{\tocsubsection}[3]{%
  \indentlabel{\@ifnotempty{#2}{\hspace*{2.3em}\makebox[2.3em][l]{%
    \ignorespaces#1 #2.\hfill}}}#3}
\renewcommand{\tocsubsubsection}[3]{%
  \indentlabel{\@ifnotempty{#2}{\hspace*{4.6em}\makebox[3em][l]{%
    \ignorespaces#1 #2.\hfill}}}#3}
\makeatother 
\begin{document}
\title{Polynomial Bound and Nonlinear Smoothing for the Benjamin-Ono Equation on the Circle}
\author{Bradley Isom, Dionyssios Mantzavinos$^*$, Seungly Oh \& Atanas Stefanov}
\begin{abstract}
For initial data in Sobolev spaces $H^s(\mathbb T)$, $\frac 12 < s \leqslant 1$, the solution to the Cauchy problem for the Benjamin-Ono equation  on the circle is shown to grow at most polynomially in time at a rate  $(1+t)^{3(s-\frac 12) + \epsilon}$, $0<\epsilon \ll 1$. 
Key to establishing this result is the discovery of a nonlinear smoothing effect for the Benjamin-Ono equation, according to which the solution to the equation satisfied by a certain gauge transform, which is widely used in the well-posedness theory of the Cauchy problem,  becomes smoother once its free solution is removed.
\end{abstract}
\date{January 19, 2020. $^*$\!\textit{Corresponding author}: mantzavinos@ku.edu.}
%
%
\thanks{B. Isom was partially supported as a graduate research assistant from NSF-DMS 1614734 and NSF-DMS 1908626. 
\\
\hspace*{3mm} A. Stefanov was partially supported by NSF-DMS 1908626.}
\maketitle
\markboth
{Polynomial Bound and Nonlinear Smoothing for the Benjamin-Ono Equation on the Circle}
{B. Isom, D. Mantzavinos, S. Oh \& A. Stefanov}
%

%
%
%
%
\section{Introduction and Results}

We consider  the Cauchy problem for the Benjamin-Ono (BO) equation on the circle
\begin{subequations}\label{bo-ivp}
\ddd{
\label{bo-eq}
&
u_t+\mathcal{H}u_{xx}=\tfrac{1}{2}\partial_x(u^2), \quad (x,t)\in\mathbb T\times\mathbb R,
\\
\label{id1}
&
u(x,0)=u_0(x) \in H^s(\mathbb T).
}
\end{subequations}
In the above initial value problem, $u(x, t)$ is a real-valued function and  $H^s(\mathbb T)$ denotes the standard $L^2$-based Sobolev space on the circle.
%
Furthermore, $\mathcal{H}$ denotes the Hilbert transform defined by
\eee{
\widehat{\mathcal{H}f}(\xi)
=
-i\, \text{sgn}(\xi)\, \widehat{f}(\xi), \quad \xi\in\mathbb Z,
}
where $\widehat f(\xi) = \mathcal F f(\xi) := \int_{x\in \mathbb T} e^{-i\xi x} f(x) dx$ is the usual Fourier transform over $\mathbb T$ and where we use the convention $\text{sgn}(0) = 0$.

The BO equation was derived in \cite{b1967, o1975} as a model for the propagation of one-dimensional long internal gravity waves in deep stratified fluids. The equation is a completely integrable system; in particular, it admits $N$-soliton solutions \cite{c1979,clp1979}, it can be expressed in the form of a Lax pair \cite{n1979,bk1979}, and it possesses an infinite number of commuting symmetries and conservation laws \cite{bk1979,ff1981}, including
\eee{
\int_{x\in \mathbb T} u dx, \quad \int_{x\in \mathbb T} u^2 dx, \quad \int_{x\in \mathbb T} \left(u\mathcal{H}u_x-\tfrac{1}{3} u^3\right) dx.
}

Without loss of generality, throughout this work we restrict our attention to solutions of the BO equation with zero mean, i.e. we assume that 
\eee{
\int_{x\in \mathbb T} u(x, t) dx =: \widehat u(0, t) = 0 \quad \forall t\in \mathbb R.
} 
This is possible thanks to the observation that  the function $v(x, t):= u(x, t)-c$ with $c = \widehat u_0(0)/2\pi$ satisfies 
$
v_t+\mathcal{H}v_{xx}
=
\partial_x(v^2)+2c v_x,
$
and hence the function $V(x,t) := v(x-2ct,t)$ satisfies the BO equation and has mean-zero initial value $V(x,0)=u_0(x)-c$. Therefore, noting that the mean $\widehat{u}(0, t)/2\pi$ of any smooth solution to the BO equation  is conserved, we deduce that $V(x,t)$ has zero mean at all times.

The Cauchy problem for the BO equation has been studied extensively in the literature. 
In the case of the line, Fokas and Ablowitz \cite{fa1983} analyzed this problem via the inverse scattering transform method under the assumption of sufficiently smooth and decaying initial data.
Iorio \cite{i1986} established local and global existence of solution for initial data in $H^s(\mathbb R)$ with $s>\frac 32$  and  $s\geqslant 2$, respectively, using energy methods (see also \cite{abfs1989}, where the continuity of the data-to-solution map is specifically addressed). 
Furthermore, Ponce \cite{p1991} proved global well-posedness for $s=\frac 32$. This result was improved by  Koch and Tzvetkov \cite{kt2003} and Kenig and Koenig \cite{kk2003}  to $s>\frac 54$ and $s>\frac 98$, respectively. 
Importantly, for initial data in $H^s(\mathbb R)$ with $s>0$, Koch and Tzvetkov \cite{kt2005} proved that the data-to-solution map of the BO initial value problem  is not uniformly continuous (previously, Molinet, Saut and Tzvetkov \cite{mst2001} had shown that the data-to-solution map fails to be $C^2$ in $H^s(\mathbb R)$  for all $s\in\mathbb R$).
This fact is due to the presence of a  derivative in the nonlinear part of the BO equation in combination with the weak smoothing effects of the linear part of the equation, and prevents one from solving the BO Cauchy problem  via a direct application of the contraction mapping principle (see also the relevant discussion in \cite{s2018}).
In \cite{t2004}, Tao bypassed this difficulty by introducing a gauge transform of Cole-Hopf type, thereby establishing global well-posedness in $H^1(\mathbb R)$.
This breakthrough idea was further employed by Burq and Planchon \cite{bp2008} and by Ionescu and Kenig \cite{ik2007}, who extended Tao's result to initial data in $H^s(\mathbb R)$ with $s>\frac 14$ and $s\geqslant 0$, respectively.
 
In the periodic setting, using Tao's gauge transform, Molinet  proved well-posedness of the Cauchy problem \eqref{bo-ivp}  in $H^s(\mathbb T)$ for $s\geqslant \frac 12$ \cite{m2007} and $s\geqslant 0$ \cite{m2008}. Furthermore, adapting the technique of \cite{kt2005} for the line, in \cite{m2007} Molinet showed that the data-to-solution map for problem \eqref{bo-ivp} is not  uniformly continuous in $H^s(\mathbb T)$ for any  $s>0$ (the corresponding result for $s<-\frac 12$ was proved by  Biagioni and Linares \cite{bl2001}). Nevertheless, Lipschitz continuity is retained in the case of mean-zero initial data (see also \cite{s2018}). 
%
Simpler proofs of the results of  \cite{m2007,m2008} along with stronger uniqueness results were provided by Molinet and Pilod in \cite{mp2012}, where an alternative proof of the result of  \cite{ik2007} was also presented. 
Finally, in the recent preprint \cite{gkt2019} G\'erard, Kappeler and Topalov obtain global well-posedness results for initial data in $H^s(\mathbb T)$ with $-\frac 12< s < 0$ (the discontinuity of the solution map for $s<-\frac 12$ had already been observed by Angulo Pava and Hakkaev \cite{aph2010}).

In order to summarize the main  results of \cite{m2007,m2008,mp2012}, we first introduce some useful notation.
\vskip 1mm
\begin{enumerate}[label=$\bullet$, leftmargin=4mm, rightmargin=0mm]
\advance\itemsep 1mm
\item
For $a,b>0$, we write $a\lesssim b$ if there exists $C>0$ such that $a\leqslant C b$. If $a\lesssim b$ and $b\lesssim a$ then we write $a\simeq b$. 
\item 
We define the Bessel potential $J_x^s$ via Fourier transform as
\sss{
\eee{
\widehat{J_x^s f}(\xi)
:=
\langle{\xi\rangle}^s\widehat{f}(\xi), \quad \langle{\cdot\rangle} := \left(1+|\cdot|^2\right)^{\frac 12}.
}
Then, for any $s\geqslant 0$ and $p\geqslant 1$, we define the Bessel potential space
\eee{
W^{s,p}(\mathbb T) := \left\{f\in L^p(\mathbb T) : \no{f}_{W^{s,p}(\mathbb T)} := \no{J_x^s f}_{L^p(\mathbb T)} <\infty \right\},
}
}
which becomes the Sobolev space $H^s(\mathbb T)$ in the special case $p=2$.
\item
For any $s, b  \in \mathbb R$, we define the Bourgain space $X_{\tau = -\omega(\xi)}^{s, b}$ by
\eee{
X_{\tau=-\omega(\xi)}^{s, b}  := \left\{f\in \mathcal D'(\mathbb T \times \mathbb R): 
\no{f}_{X^{s, b}} := \no{\left\langle \xi \right\rangle^s \left\langle \tau+\omega(\xi)\right\rangle^b \widetilde f(\xi, \tau)}_{L^2(\mathbb Z_\xi\times \mathbb R_\tau)} <\infty \right\},
}
where $\widetilde f(\xi, \tau)$  denotes the Fourier transform of $f(x, t)$ with respect to both $x$ and $t$.
We denote by $X_{\tau=-\omega(\xi), T}^{s, b}$ the restriction of $X_{\tau=-\omega(\xi)}^{s, b}$ on $\mathbb T \times [0, T]$.
Furthermore, for convenience of notation, hereafter we shall write $X^{s,b}_{\tau=-\xi^2} =: X^{s,b}$ and $X^{s,b}_{\tau=\xi^2} =: \bar X^{s,b}$ and, analogously, $X_T^{s,b}$ and $\bar X_T^{s,b}$ for the  respective restrictions of these spaces on $\mathbb T \times [0, T]$.
\item
The operators $\Pi^0$, $\Pi^+$ and $\Pi^-$ denote the projections onto the zero, positive and negative Fourier modes, respectively: 
\eee{
\widehat{\Pi^0(f)}
:=
\frac{1}{2\pi} \widehat{f}(0),
\quad
\widehat{\Pi^{\pm}(f)}(\xi)
:=
\chi^{\pm}(\xi)\widehat{f}(\xi),
}
where $\chi^\pm(\xi)$ are the characteristic functions for $\xi\gtrless 0$. It is straightforward to see that
$
\overline{\Pi^{\pm}(u)}=\Pi^{\mp}(\overline{u}).
$
\item
For $k\geqslant 1$, we define the Littlewood-Paley-type projection operator $P_k$ by 
\eee{
\widehat{P_k(f)}(\xi)
:=
\chi_{\left\{2^{k-1} \leqslant |\xi|< 2^k\right\}}\widehat{f}(\xi),
}
where $\chi_A$ is the characteristic function of the set $A$.  We will often denote $P_k(f)$ simply by $f_k$. By this definition, it follows that
\eee{
2\pi \, \widehat{\Pi^0(f)}(\xi) +\sum_{k=1}^{\infty}\widehat{f_k}(\xi)
=
\widehat{f}(\xi), \quad \xi\in\mathbb Z.}
\item
Following \cite{mp2012}, we introduce the gauge transform for the periodic BO equation as
\eee{\label{w-def}
w := \p_x  \Pi^+  (e^{-iF/2}),
}
where $F = \p_x^{-1} u$ is the primitive of the solution $u$ of problem  \eqref{bo-ivp} such that
\eee{
\widehat{F}(\xi, t)
:=
\def\arraystretch{1.35}
\left\{ 
\begin{array}{ll} 
0, & \xi=0, 
\\ 
\dfrac{1}{i\xi} \, \widehat{u}(\xi, t), & \xi\in\mathbb Z\setminus\{0\}. 
\end{array} 
\right.
}
Noting that $F$ has zero mean and is $2\pi$-periodic, it is straightforward to see that it satisfies the equation
\eee{\label{F-pde}
F_t + \mathcal{H}F_{xx} - \frac 12 F_x^2
= -\frac 12 \widehat{\Pi^0(F_x^2)}(t).
}
In turn, noting also that for any mean-zero function $f$ we have 
$\mathcal{H}f = -if+2i \Pi^-(f)$, we infer that $w$ satisfies the initial value problem
\sss{\label{w-ivp}
\ddd{
&w_t-iw_{xx} = -\partial_x\Pi^+\big((\partial_x^{-1}w) \Pi^-(u_x)\big)
+ \tfrac i4 \widehat{\Pi^0(u^2)}w, \quad &&(x,t)\in\mathbb T\times\mathbb R,
\label{w-eq}
\\
&w(x, 0) = \p_x  \Pi^+ (e^{-i\p_x^{-1} u_0(x)/2}) =:  w_0(x), && x\in \mathbb T.
\label{w-ic}
}
}
\end{enumerate}

With the above definitions at hand, the main well-posedness results of \cite{m2007,m2008,mp2012}  can be summarized as follows (see, in particular, Theorem 7.1 in \cite{mp2012}).
\begin{theorem}[\b{Well-posedness on the circle} \cite{m2007,m2008,mp2012}]
\label{bo-wp-t}
Suppose $u_0\in H^s(\mathbb T)$ with $0\leqslant s\leqslant 1$. Then,  the initial value problem \eqref{bo-ivp} for the BO equation on the circle admits a solution
$$
u\in C\big([0,T];H^s(\mathbb T)\big) \cap L^4\big([0,T];W^{s,4}(\mathbb T)\big) \cap X^{s-1,1}_{\tau=-|\xi|\xi, T}
$$
where $T = T\big(\no{u_0}_{L^2(\mathbb T)}\big) \simeq \min\big\{\no{u_0}_{L^2(\mathbb T)}^{-4},1\big\} >0$.
Moreover, the initial value problem \eqref{w-ivp} for the function $w$, which is defined in terms of $u$ via the gauge transform \eqref{w-def}, admits a solution  
$$
w\in C([0,T];H^s(\mathbb T))\cap X^{s,\frac 12}_T
$$
in the distributional as well as in the Duhamel sense.
In particular, we have the estimates
\ddd{
&\max\Big\{\no{u}_{C([0,T];H^s(\mathbb T))}, \no{u}_{L^4([0,T];W^{s,4}(\mathbb{T}))}, \no{w}_{X^{s,\frac 12}_T}\Big\}
\lesssim
\max\left\{\no{u_0}_{L^2(\mathbb T)}^{2s}, 1 \right\} \no{u_0}_{H^s(\mathbb T)},
\nn\\
&\no{u}_{X_{\tau=-|\xi|\xi, T}^{s-1, 1}}
\lesssim
 \left( \no{u_0}_{H^s(\mathbb T)} 
+  \no{u_0}_{H^s(\mathbb T)}^2 \right).
\nn
}
\end{theorem}

\begin{remark}[\b{Global well-posedness}]
Thanks to the conservation of the $L^2$-norm, the solution of Theorem \ref{bo-wp-t} is in fact a global solution \cite{m2007,m2008} which is unique within the class of limits of smooth solutions of problem \eqref{bo-ivp}.
\end{remark}

The scope of the present work extends beyond the fundamental question of well-posedness that was addressed in \cite{m2007,m2008,mp2012}.
In particular, we revisit the Cauchy problem \eqref{bo-ivp} for the BO equation on the circle and obtain an explicit growth bound of polynomial type for the solution guaranteed by Theorem \ref{bo-wp-t}. 
Crucial for proving this bound is a nonlinear smoothing effect that we establish for the BO equation, according to which the nonlinear component of the solution of the equation emanating from the gauge transform is smoother than the component corresponding to the initial datum.
More precisely, we shall show the following.
\begin{theorem}[\b{Nonlinear smoothing}]
\label{nl-smooth-t}
Suppose $\frac 16 < s \leqslant 1$, $0<a< \min\left\{s-\frac 16,\frac{1}{3}\right\}$ and $K:=\frac{1}{8\pi} \no{u_0}^2_{L^2(\mathbb T)}$, and let $u$ and $w$ be the solutions of the Cauchy problems \eqref{bo-ivp} and \eqref{w-ivp} established by Theorem \ref{bo-wp-t}. Then, 
$
e^{-iKt} w(x, t)- e^{it\partial_x^2}w_0(x)
\in C([0,T];H^{s+a}(\mathbb T))
$
with the estimate
$$
\no{e^{-iKt} w - e^{it\partial_x^2}w_0}_{C([0,T];H^{s+a}(\mathbb T))}
\leqslant
C\big(\no{u_0}_{H^{\min\{ s, \frac 12\}}(\mathbb T)}\big)  \no{u_0}_{H^s(\mathbb T)},
$$
where $e^{it\partial_x^2}$ is the semigroup associated with the linear Schr\"odinger equation.
\end{theorem}

The nonlinear smoothing effect of Theorem \ref{nl-smooth-t} provides the basis for proving the following polynomial growth bound for the solution of the BO initial value problem \eqref{bo-ivp}.

\begin{theorem}[\b{Polynomial bound}]
\label{poly-bound-t}
Suppose $\frac 12 < s \leqslant 1$.  Then, for any $0<\epsilon \ll 1$, the solution $u$ of the BO Cauchy problem \eqref{bo-ivp} established by Theorem \ref{bo-wp-t} satisfies 
\eee{
\n{u(t)}{H^s(\mathbb T)}  \leqslant C\big(\epsilon,s, \n{u_0}{H^s(\mathbb T)}\big) \lan{t}^{3 (s - \frac 12) + \epsilon}, \quad t\in \mathbb R.
}
 \end{theorem}

The connection between nonlinear smoothing and polynomial bounds for Hamiltonian equations was first established by Bourgain \cite{b1996,b1997}, who employed Fourier truncation operators in conjunction with smoothing estimates to obtain the following local-in-time inequality for solutions of various dispersive PDEs:
\begin{equation}\label{eq:bour}
\n{ u(t + \delta) }{H^s} \leqslant \n{u(t)}{H^s} + C\n{u(t)}{H^s}^{1-\delta}
\end{equation}
 for some $\delta \in (0,1)$.  Local time iterations using the above inequality resulted in the polynomial growth bound $\n{u(t)}{H^s} \lesssim \left\langle t \right\rangle^{1/\delta}$.  Staffilani \cite{s1997b, s1997a} used further multilinear smoothing estimates to obtain \eqref{eq:bour} which led to polynomial bounds of high-Sobolev norms $s>1$ for Korteweg-de Vries (KdV) and nonlinear Schr\"odinger (NLS) equations.  Colliander, Keel, Staffilani, Takaoka and Tao \cite{ckstt2003} developed a new method using modified energy called the ``upside-down $I$-method'' to produce polynomial bounds in low-Sobolev norms $s\in (0,1)$ for the NLS equation.  Sohinger \cite{s2011a, s2011b} further developed the upside-down I-method to obtain polynomial bounds for high Sobolev norms for NLS.  We also refer the reader to \cite{cko2012} and the references therein for further developments. 

More recently, uniform-in-time bounds have been established for a number of completely integrable dispersive equations using inverse scattering techniques. In particular, Killip, Visan and Zhang \cite{kvz2018} showed that the $H^s$-norm of solutions to the KdV and NLS equations is uniformly bounded in time for $-1 \leqslant s < 1$ and $-\frac 12 < s < 1$, respectively, both on the line and on the circle.  For the BO equation,  Talbut \cite{t2019} proved an analogous bound  in $H^s$ for $-\f{1}{2} < s < 0$.  However, no bound is available for $s>0$ since the  technique used in \cite{t2019}, which is similar to that of \cite{kvz2018},  becomes rather convoluted for higher values of $s$.  On the other hand, Koch and Tataru \cite{kt2018} showed that there exists a conserved energy equivalent to the $H^s$-norm for $s>-\f{1}{2}$ in the case of the NLS and mKdV equations and for $s\geqslant -1$ in the case of the KdV equation.

Nonlinear smoothing properties analogous to the one of Theorem \ref{nl-smooth-t} have been previously established for several important dispersive equations. Indicatively, we mention the work of Erd\u ogan and Tzirakis \cite{et2013a} on the periodic Korteweg-de Vries (KdV) equation, as well as their  works on the derivative nonlinear Schr\"odinger equation on the line \cite{egt2018}, the fractional Schr\"odinger equation on $\mathbb T$ and $\mathbb R$ \cite{egt2019}, and the Zakharov system on the torus \cite{et2013b}.
The main technique used in the proof of these results is known as the normal form method and was first introduced by Shatah \cite{s1985} in the context of the Klein-Gordon equation with a quadratic nonlinearity. This method was further developed recently by Germain, Masmoudi and Shatah  for two-dimensional quadratic Schr\"odinger equations \cite{gms2009a} and the gravity water waves equation \cite{gms2009b}, as well as by  Babin, Ilyin and Titi  for the periodic KdV equation \cite{bit2011}. The technique used in the latter work is known as  differentiation by parts. An alternative formulation of the normal form method which involves a multilinear pseudo-differential operator in place of differentiation by parts was provided in \cite{os2012, o2013}. In our work, thanks to the properties of the gauge transform,  the normal form machinery is not required for proving the nonlinear smoothing result of Theorem \ref{nl-smooth-t}.

\vskip 2mm
\noindent
\textbf{Structure of the paper.} 
In Section \ref{nl-est-s} we establish a bilinear estimate which is crucial for showing the nonlinear smoothing effect of Theorem \ref{nl-smooth-t}. The proof of this theorem is then provided in Section \ref{nl-smooth-s}. Finally, the polynomial growth bound of Theorem \ref{poly-bound-t} is established in Section \ref{poly-bound-s}.

\section{Bilinear Estimate}\label{nl-est-s}

The following bilinear estimate plays a key role in the proof of the nonlinear smoothing effect of Theorem \ref{nl-smooth-t}.
\begin{proposition}[\b{Bilinear estimate}]
\label{main-est-p}
Let $V\in X^{0,\frac 12}$ and $U\in L^\infty(\mathbb R; L^2(\mathbb T))\cap \bar X^{0,1}$ with $U$ compactly supported in $[-T,T]$ for some $T>0$. Then, for all $\delta > 0$, $m\in \mathbb N$, $k\in \mathbb N$ and $0<j\leqslant k$, we have
\ddd{\label{smooth-est}
&\quad
\no{P_j\Pi^+ (V_k \, \Pi^-(U_m))}_{X^{0,-\frac 12-\delta}} 
\nn\\
&\lesssim
2^{\left(\frac 16+\delta\right)k-\frac{m+j}{2}}
\no{V_k}_{X^{0,\frac 12}}
\Big(\no{\Pi^-(U_m)}_{L^\infty(\mathbb R; L^2(\mathbb T))}
+2^{-\frac{m+j}{2}} \no{\Pi^-(U_m)}_{\bar X^{0,1}}\Big),
}
where the implicit constant depends on $T$.
\end{proposition}

Indeed, via complex interpolation it can be shown that
$\big(X^{0, -\frac 12 -\delta}, X^{0, 0}\big)_\theta
=
X^{0, -\frac 12 + \delta}$ with $\theta :=\frac{1/2-\delta}{1/2+\delta}$.
Therefore, interpolating between estimate \eqref{smooth-est} and the estimate
$$
\no{P_j\Pi^+(V_k \Pi^-(U_m))}_{L^2(\mathbb T \times \mathbb R)}
\lesssim
\no{V_k}_{X^{0,\frac 12}} \no{\Pi^-(U_m)}_{\bar X^{0,1}},
$$
which follows from the generalized H\"older inequality and the embedding \cite{b1993}
$X_{\tau=\pm \xi^2}^{0,\frac{3}{8}+\delta}\hookrightarrow L^4(\mathbb T\times\mathbb R)$, $\delta>0$,  we obtain
\ddd{\label{interp-c}
\no{P_j\Pi^+(V_k \Pi^-(U_m))}_{X^{0,-\frac 12+\delta}}
&\lesssim
\left[
2^{-\frac{m+j}{2}+k(\frac 16+\delta)} 
\left(\no{\Pi^-(U_m)}_{L^\infty(\mathbb R; L^2(\mathbb T))}
+
2^{-\frac{m+j}{2}} 
\no{\Pi^-(U_m)}_{\bar X^{0,1}}
\right)
\right]^{\theta}
\nn\\
&\quad
\cdot
\no{\Pi^-(U_m)}_{\bar X^{0,1}}^{1-\theta} 
\no{V_k}_{X^{0,\frac 12}}.
}
This estimate is a main ingredient in the proof of Theorem \ref{nl-smooth-t}  which is provided in Section \ref{nl-smooth-s}. In the remaining of the current section, we prove Proposition \ref{main-est-p}.

\begin{Proof}[Proof of Proposition \ref{main-est-p}]
Observe that $\varphi\in X^{s,b}$ implies $\bar{\varphi}\in\bar X^{s,b}$. 
By the dual formulation of the Bourgain norm along with Plancherel's theorem, we have
\ddd{
&\quad
\no{P_j\Pi^+(V_k \, \Pi^-(U_m))}_{X^{0,-\frac 12-\delta}}
=
\sup_{\no{\varphi}_{X^{0,\frac 12+\delta}}=1}
\left|
\int_{x \in \mathbb T} \int_{t\in \mathbb R}
P_j\Pi^+(V_k \, \Pi^-(U_m)) \cdot \overline{\varphi} (x, t)  dt dx
\right|
\nn\\
&\simeq
\sup_{\no{\varphi}_{X^{0,\frac 12+\delta}}=1}
\Bigg|
\sum_{\xi_1 \in \mathbb Z} \sum_{\xi_2 \in \mathbb Z} \int_{\tau_1\in \mathbb R} \int_{\tau_2 \in \mathbb R}
\widetilde{V_k}(\xi_1,\tau_1) 
\widetilde{\Pi^-(U_m)}(\xi_2,\tau_2) 
\widetilde{\Pi^-(\overline \varphi_j)}(-(\xi_1+\xi_2),-(\tau_1+\tau_2)) d\tau_2 d\tau_1
\Bigg|.
\label{dual-norm-prop}
}
Next, we let 
$L_1=\left|\tau_1+\xi_1^2\right|$, $L_2=\left|\tau_2-\xi_2^2\right|$, $L_3=\big|-\left(\tau_1+\tau_2\right)-\left(-\left(\xi_1+\xi_2\right)\right)^2\big|$ and observe that, since 
$
-\left(\tau_1+\tau_2\right)-\left(\xi_1+\xi_2\right)^2
=
-(\tau_1+\xi_1^2)-(\tau_2-\xi_2^2)-2\xi_2(\xi_1+\xi_2)
$
and $2^{m-1}\leqslant |\xi_2| < 2^m$, $2^{j-1}\leqslant \left|\xi_1+\xi_2\right| < 2^j$, we have
\eee{\label{max-ineq}
\max \left\{L_1,L_2,L_3\right\} \geqslant \tfrac 16 \, 2^{m+j}.
}
Then, writing $1 = \chi_{A_1} +  \chi_{A_1^c} \chi_{A_3} + \chi_{A_1^c}\chi_{A_2}  \chi_{A_3^c}$ with 
$A_i := \big\{L_i\geqslant \frac 16 \, 2^{m+j}\big\}$, $i=1,2,3$, we have
\eee{\label{I123}
\bigg|
\sum_{\xi_1 \in \mathbb Z} \sum_{\xi_2 \in \mathbb Z} \int_{\tau_1\in \mathbb R} \int_{\tau_2 \in \mathbb R}
\widetilde{V_k}(\xi_1,\tau_1) 
\widetilde{\Pi^-(U_m)}(\xi_2,\tau_2) 
\widetilde{\Pi^-(\overline \varphi_j)}(-(\xi_1+\xi_2),-(\tau_1+\tau_2)) d\tau_2 d\tau_1
\bigg|
\leqslant
I_1 + I_2 + I_3
}
where
\ddd{
&I_1 =
\bigg|
\sum_{\xi_1 \in \mathbb Z} \sum_{\xi_2 \in \mathbb Z} \int_{\tau_1\in \mathbb R} \int_{\tau_2 \in \mathbb R}
\chi_{A_1} 
\widetilde{V_k}(\xi_1,\tau_1) 
\widetilde{\Pi^-(U_m)}(\xi_2,\tau_2) 
\widetilde{\Pi^-(\overline \varphi_j)}(-(\xi_1+\xi_2),-(\tau_1+\tau_2)) d\tau_2 d\tau_1
\bigg|
\nn\\
&I_2 =
\bigg|
\sum_{\xi_1 \in \mathbb Z} \sum_{\xi_2 \in \mathbb Z} \int_{\tau_1\in \mathbb R} \int_{\tau_2 \in \mathbb R}
\chi_{A_1^c} 
\widetilde{V_k}(\xi_1,\tau_1) 
\widetilde{\Pi^-(U_m)}(\xi_2,\tau_2) 
\chi_{A_3}
\widetilde{\Pi^-(\overline \varphi_j)}(-(\xi_1+\xi_2),-(\tau_1+\tau_2)) d\tau_2 d\tau_1
\bigg|
\nn\\
&I_3 =
\bigg|
\sum_{\xi_1 \in \mathbb Z} \sum_{\xi_2 \in \mathbb Z} \int_{\tau_1\in \mathbb R} \int_{\tau_2 \in \mathbb R}
\chi_{A_1^c} 
\widetilde{V_k}(\xi_1,\tau_1) 
\chi_{A_2}
\widetilde{\Pi^-(U_m)}(\xi_2,\tau_2) 
\chi_{A_3^c}
\widetilde{\Pi^-(\overline \varphi_j)}(-(\xi_1+\xi_2),-(\tau_1+\tau_2)) d\tau_2 d\tau_1
\bigg|.
\nn
}

We begin with the estimation of $I_1$.
Define $f^{A_i}$ via its Fourier transform as $\widetilde{f^{A_i}} := \chi_{A_i}\widetilde{f}$. 
Then, Plancherel's theorem followed by the Cauchy-Schwarz inequality yield
\ddd{
I_1 
&\simeq
\left|
\int_{x \in \mathbb T} \int_{t\in \mathbb R} 
V_k^{A_1}(x, t) 
\cdot
\Pi^-(U_m)(x,t) 
\cdot
\Pi^-(\overline \varphi_j)(x, t) dt dx
\right|
\nn\\
&\leqslant
\big\|V_k^{A_1}\big\|_{L^2(\mathbb T \times \mathbb R)} 
\no{\Pi^-(U_m) \cdot \Pi^-(\overline \varphi_j)}_{L^2(\mathbb T \times \mathbb R)}.
\label{prop-temp0}
}
For the first factor in \eqref{prop-temp0}, recalling the definition of $A_1$ we proceed as follows:
\ddd{
\big\|V_k^{A_1}\big\|_{L^2(\mathbb T \times \mathbb R)}
&\simeq
\bigg(
\sum_{\xi \in \mathbb Z}\int_{\tau\in \mathbb R}
\frac{|\tau+\xi^2|}{\left|\tau+\xi^2\right|}
\, \big|\widetilde{V_k^{A_1}}(\xi,\tau)\big|^2 d\tau
\bigg)^{\frac 12}
\nn\\
&\leqslant
\bigg(
\sum_{\xi \in \mathbb Z}\int_{\tau\in \mathbb R}
\frac{\big(1+|\tau+\xi^2|^2\big)^{\frac 12}}{\frac 16\, 2^{m+j}}
\, \big|\widetilde{V_k^{A_1}}(\xi,\tau)\big|^2 d\tau
\bigg)^{\frac 12}
\simeq
2^{-\frac{m+j}{2}}
\no{V_k}_{X^{0,\frac 12}}.
\label{prop-temp1}
}
For the second factor in \eqref{prop-temp0}, recalling that $U$ is supported inside $[-T, T]$ and applying the generalized H\"older inequality, we have
\ddd{
\no{\Pi^-(U_m)\cdot\Pi^-(\overline \varphi_j)}_{L^2(\mathbb T \times \mathbb R)}
&\leqslant
\no{\Pi^-(U_m)}_{L^\infty([-T,T]; L^2(\mathbb T))}
\no{\Pi^-(\overline \varphi_j)}_{L^2([-T,T];L^\infty(\mathbb T))}
\nn\\
&\lesssim
T^{\frac 13}
\no{\Pi^-(U_m)}_{L^\infty([-T,T]; L^2(\mathbb T))}
\no{\Pi^-(\overline \varphi_j)}_{L^6([-T,T]; L^\infty(\mathbb T))}.
\nn
}
Moreover,   the Sobolev embedding $W^{\sigma,p}(\mathbb T)\hookrightarrow L^\infty(\mathbb T)$, $1\leqslant p\leqslant\infty$, $\sigma>\frac{1}{p}$, for $p=6$  yields
$$
\no{\Pi^-(\overline \varphi_j)}_{L^6([-T,T]; L^\infty(\mathbb T))}
\lesssim
\big\|J_x^\sigma \Pi^-(\overline \varphi_j)\big\|_{L^6([-T,T]; L^6(\mathbb T))}, \quad \sigma>\tfrac 16,
$$
while the embedding $X_{\tau=\pm \xi^2}^{\varepsilon,\frac 12+\delta} \hookrightarrow L^6(\mathbb T \times\mathbb R)$,  $\varepsilon, \delta>0$   \cite{b1993} further implies
\eee{
\no{\Pi^-(\overline \varphi_j)}_{L^6([-T,T]; L^\infty(\mathbb T))}
\lesssim
\big\|J_x^\sigma \Pi^-(\overline \varphi_j)\big\|_{\bar X^{\varepsilon,\frac 12+\delta}}
\lesssim
2^{j(\sigma+\varepsilon)}
\no{\overline{\varphi}}_{\bar X^{0,\frac 12+\delta}}.
\nn
}
In turn, we find
\eee{\label{prop-temp2}
\no{\Pi^-(U_m)\cdot\Pi^-(\overline \varphi_j)}_{L^2(\mathbb T \times \mathbb R)}
\lesssim
T^{\frac 13}
\no{\Pi^-(U_m)}_{L^\infty([-T,T]; L^2(\mathbb T))}
2^{j(\sigma+\varepsilon)}
\no{\overline{\varphi}}_{\bar X^{0,\frac 12+\delta}}.
}
Hence,  setting $\sigma+\ve = \frac 16 + \delta$ with $\delta>\ve$ and then combining \eqref{prop-temp2} with \eqref{prop-temp0} and \eqref{prop-temp1}, we deduce
\eee{\label{I1-est}
I_1
\lesssim
2^{-\frac{m+j}{2} + j\left(\frac 16+\delta\right)}
\no{V_k}_{X^{0,\frac 12}}
\no{\Pi^-(U_m)}_{L^\infty([-T,T]; L^2(\mathbb T))}
\no{\overline{\varphi}}_{\bar X^{0,\frac 12+\delta}}.
}

We continue with the estimation of $I_2$. As with $I_1$, we employ Plancherel's theorem and the Cauchy-Schwarz inequality  to infer
\eee{\label{I2-temp1}
I_2
\leqslant
\no{(\Pi^-(\overline \varphi_j))^{A_3}}_{L^2(\mathbb T \times \mathbb R)}
\big\|V_k^{A_1^c}\, \Pi^-(U_m)\big\|_{L^2(\mathbb T \times \mathbb R)}.
}
Then, similarly to \eqref{prop-temp1} we have
$$
\no{(\Pi^-(\overline \varphi_j))^{A_3}}_{L^2(\mathbb T \times \mathbb R)}
\lesssim
2^{-(m+j)(\frac 12+\delta)}
\no{\overline{\varphi}}_{\bar X^{0,\frac 12+\delta}}.
$$
Moreover, treating the second factor in \eqref{I2-temp1} similarly to the corresponding term in $I_1$, we find
$$
\big\|V_k^{A_1^c} \, \Pi^-(U_m)\big\|_{L^2(\mathbb T \times \mathbb R)}
\lesssim
2^{k(\frac 16 + \delta)}
\no{\Pi^-(U_m)}_{L^\infty([-T,T]; L^2(\mathbb T))}
\big\|V_k^{A_1^c}\big\|_{X^{0,\frac 12+\delta}}.
$$
Hence, observing that
$
\|V_k^{A_1^c}\|_{X^{0,\frac 12+\delta}} \lesssim
2^{\delta(m+j)} \no{V_k}_{X^{0,\frac 12}}
$
by the definition of $A_1$, we conclude that
\eee{\label{I2-est}
I_2
\lesssim
2^{-\frac{m+j}{2} + k(\frac 16+\delta)}
\no{V_k}_{X^{0,\frac 12}} 
\no{\Pi^-(U_m)}_{L^\infty([-T,T]; L^2(\mathbb T))} 
\no{\overline{\varphi}}_{\bar X^{0,\frac 12+\delta}}.
}

Finally, similarly to $I_1$ and $I_2$, for $I_3$ we have
\eee{\label{I3-temp1}
I_3
\leqslant
\no{(\Pi^-(U_m))^{A_2}}_{L^2(\mathbb T \times \mathbb R)}
\big\|V_k^{A_1^c} \big(\Pi^-(\overline \varphi_j)\big)^{A_3^c}\big\|_{L^2(\mathbb T \times \mathbb R)}.
}
For the first factor in \eqref{I3-temp1}, we proceed as with \eqref{prop-temp1} to obtain
\eee{\label{I3-temp3}
\no{(\Pi^-(U_m))^{A_2}}_{L^2(\mathbb T \times \mathbb R)}
\lesssim
2^{-(m+j)} 
\no{\Pi^-(U_m)}_{\bar X^{0,1}}.
}
Moreover, for the second factor in \eqref{I3-temp1}, we use the generalized H\"older inequality as well as the embedding $X_{\tau=\pm \xi^2}^{0,\frac{3}{8}+\delta}\hookrightarrow L^4(\mathbb T\times\mathbb R)$ to find
\eee{
\big\|V_k^{A_1^c} \big(\Pi^-(\overline \varphi_j)\big)^{A_3^c}\big\|_{L^2(\mathbb T \times \mathbb R)}
\leqslant
\big\|V_k^{A_1^c}\big\|_{L^4(\mathbb T \times \mathbb R)} 
\no{(\Pi^-(\overline \varphi_j))^{A_3^c}}_{L^4(\mathbb T \times \mathbb R)}
\lesssim
\no{V_k}_{X^{0,\frac 12}} 
\no{\overline{\varphi}}_{\bar X^{0,\frac 12+\delta}}.
\label{I3-temp2}
}
Therefore, combining \eqref{I3-temp2} and \eqref{I3-temp3} into \eqref{I3-temp1}, we deduce 
\eee{\label{I3-est}
I_3
\lesssim
2^{-(m+j)} 
\no{V_k}_{X^{0,\frac 12}} 
\no{\overline{\varphi}}_{\bar X^{0,\frac 12+\delta}} 
\no{\Pi^-(U_m)}_{\bar X^{0,1}}.
}

Overall,  the three estimates \eqref{I1-est}, \eqref{I2-est} and \eqref{I3-est} together with the decomposition \eqref{I123} and the dual formulation \eqref{dual-norm-prop} imply the desired estimate \eqref{smooth-est}.
\end{Proof}

\section{Nonlinear Smoothing: Proof of Theorem \ref{nl-smooth-t}}\label{nl-smooth-s}

We begin by noting that the existence of  the solution $u$ of Theorem \ref{bo-wp-t} for the BO Cauchy problem \eqref{bo-ivp} on $\mathbb T \times [0, T]$  is proved by first taking initial data $u_0\in H^s(\mathbb T)$ with small $L^2$-norm and constructing $u$ as the strong  limit of a sequence of smooth solutions $u_n \in C([0,1];H^s(\mathbb T))\cap L^4([0,1];W^{s,4}(\mathbb T))\cap X^{s-1,1}_{\tau=-|\xi|\xi, 1}$. 
Also, in \cite{mp2012} it is shown that  the sequence of gauge transforms $w_n := \partial_x\Pi^+(e^{-iF_n/2})$ corresponding to  $u_n = \p_x F_n$ converges to some $w$ in $C([0,1];H^s(\mathbb T)) \cap X^{s,1/2}_1$. Furthermore, due to the strong convergence of   $u_n$  in $C([0,1];H^s(\mathbb T))$ it follows  from the mean value theorem that $w_n$ converges to $\partial_x\Pi^+(e^{-iF/2})$ in $C([0,1];L^2(\mathbb T))$, and hence $w=\partial_x\Pi^+(e^{-iF/2})$.
In turn, it follows that $v_n:= e^{-iKt} w_n$ converges to 
\eee{\label{v-def}
v(x, t) := e^{-iKt}w(x, t)
}
in $C([0,1];H^s(\mathbb T)) \cap X^{s,1/2}_1$. 
Then, using the smoothness of $v_n$ together with standard estimates (e.g. estimate \eqref{tao1-l}) and Proposition \ref{main-est-p}, it follows that $v$ satisfies the Duhamel equation
\eee{\label{duhamel}
v(x, t)
=
\eta(t) e^{it\partial_x^2} v_0(x)
-
\eta(t)
\int_{t'=0}^t e^{i(t-t')\partial_x^2}
\partial_x\Pi^+(\partial_x^{-1} v\cdot \partial_x\Pi^-(u))(x, t')
dt', \quad t\in [0, 1],
}
where $\eta\in C_0^{\infty}(\mathbb R)$ is supported inside $[-2,2]$ with $\eta\equiv 1$ on $[-1,1]$ and $0\leqslant \eta \leqslant 1$ for all $t\in \mathbb R$.
In addition, observe that if $u$ solves \eqref{bo-ivp} then so does $\lambda u(\lambda x, \lambda^2 t)$. Exploiting this scaling with $\lambda = 1/T^2$ and the fact that all previous convergences hold in spaces where the spatial period is assumed to be $\lambda\geqslant 1$ \cite{mp2012},  the small  $L^2$-norm assumption on $u_0$ can be dropped  and the lifespan of the solution can be extended to the lifespan $T \simeq \min \big\{\no{u_0}_{L^2(\mathbb T)}^{-4}, 1\big\}$ of Theorem \ref{bo-wp-t}. 
Therefore,  $v$ satisfies the Duhamel equation \eqref{duhamel} on $[0,T]$, i.e.
\eee{\label{duhamel-T}
v(x, t)
=
\eta_T(t)  e^{it\partial_x^2} v_0(x)
-
\eta_T(t)
\int_{t'=0}^t e^{i(t-t')\partial_x^2}\, \partial_x\Pi^+(\partial_x^{-1}v\cdot\Pi^-(u_x))(x, t') dt', \quad t\in [0, T],
}
where $\eta_T(t) := \eta(t/T)$. 

Combining the representation \eqref{duhamel-T} with the embedding 
$X_T^{s,b} \hookrightarrow C([0, T];H^s(\mathbb T))$, $s\in \mathbb R$, $b>\tfrac 12$,
(see, for example, Corollary 2.10 in \cite{t2006}) we obtain
\eee{\label{nls-temp1}
\no{e^{-iKt} w - e^{it\partial_x^2}w_0}_{C([0,T];H^{s+a}(\mathbb T))}
\lesssim
\no{\eta_T \int_{t'=0}^t e^{i(t-t')\partial_x^2}\, \partial_x\Pi^+(\partial_x^{-1}v\cdot\Pi^-(u_x)) dt'}_{X^{s+a,\frac 12+\delta}_T}.
}
In order to estimate the right-hand side of the above inequality, we first need to define appropriate extensions of the functions $v$ and $u$ with respect to $t$ outside the interval $[0,T]$. 
For $v$, we choose an extension $v^*\in X^{s,\frac 12}$ such that
\eee{\label{v-ext-def}
\no{v^*}_{X^{s,\frac 12}}
\leqslant
2\no{v}_{X^{s,\frac 12}_T},
}
which exists for all $s\in \mathbb R$ by the definition of $X_T^{s, b}$ as a restriction of $X^{s, b}$.
For $u$, we use a less trivial extension which is similar to the one in \cite{mpv2019} and is defined as follows.

\begin{lemma}[\b{Extension of $u$ outside $[0, T]$}]
\label{u-ext-def-l}
Given $u\in C([0,T];H^s(\mathbb T))\cap X_{\tau=-|\xi|\xi, T}^{s-1,1}$, let 
\eee{\label{u-ext-def}
u^*(t) := S(t)\eta_T(t) S(-\mu_T(t)) u(\mu_T(t)),
}
where $S(\cdot)$ is the free group associated with the linear component of the BO equation, whose action is defined by $\widehat{S(t)f}(\xi) := e^{-i |\xi| \xi t} \widehat f(\xi)$,  and
\eee{
\mu_T(t)=
\left\{
\begin{array}{ll}
t, & t\in [0,T],
\\
2T-t, & t\in [T,2T],
\\
0, & t\notin [0,2T].
\end{array}
\right.
\nn
}
If there exists a smooth approximating sequence $u_n$ for $u$  in $C([0,T];H^s(\mathbb T))\cap X^{s-1,1}_{\tau=-|\xi|\xi,T}$, then
\ddd{
&\no{u^*}_{L^\infty(\mathbb R; H^s(\mathbb T))}
\lesssim
\no{u}_{C([0,T];H^s(\mathbb T))},
\label{ext-est1}
\\
&\no{u^*}_{X_{\tau=-|\xi|\xi}^{s-1,1}}
\lesssim
\no{u}_{X_{\tau=-|\xi|\xi, T}^{s-1,1}}
+
\no{u}_{C([0,T];H^s(\mathbb T))},
\label{ext-est2}
}
where the implicit constants depend on $T$.
\end{lemma}

\begin{Proof}[Proof of Lemma \ref{u-ext-def-l}]
For inequality \eqref{ext-est1}, we simply note that
$$
\no{u^*}_{L^\infty(\mathbb R; H^s(\mathbb T))}
\lesssim
\no{u(\mu_T)}_{C([-2T,2T];H^s(\mathbb T))}
=
\no{u}_{C([0,T];H^s(\mathbb T))}.
$$

For inequality \eqref{ext-est2}, we let $u_n$ be an approximating sequence for $u$ in $C([0,T];H^s(\mathbb T))\cap X^{s-1,1}_{\tau=-|\xi|\xi,T}$ and denote by $u_n^*$ its extension defined analogously to \eqref{u-ext-def}. By the definition of the Bourgain norm and  the properties of $\eta_T$ and $\mu_T$, we find 
\ddd{
\no{u_n^*}_{X^{s-1,1}_{\tau=-|\xi|\xi}} 
&\lesssim
\no{S(-\mu_T) u_n(\mu_T)}_{L^2([-2T,2T];H^{s-1}(\mathbb T))}
+
\no{\partial_t(S(-\mu_T) u_n(\mu_T))}_{L^2([-2T,2T];H^{s-1}(\mathbb T))}
\nn\\
&\lesssim
\no{u_n(0)}_{H^{s-1}(\mathbb T)}
+
\no{u_n}_{L^2([0,T];H^{s-1}(\mathbb T))}
+
\no{\partial_t(S(-\mu_T) u_n(\mu_T))}_{L^2([-2T,0];H^{s-1}(\mathbb T))}
\nn\\
&\quad+
\no{\partial_t(S(-\mu_T) u_n(\mu_T))}_{L^2([0,T];H^{s-1}(\mathbb T))}
+
\no{\partial_t(S(-\mu_T) u_n(\mu_T))}_{L^2([T,2T];H^{s-1}(\mathbb T))},
\nn
}
where the implicit constant in the second inequality depends on $T$. Since $u_n$ is smooth, we directly compute
$$
\partial_t(S(-\mu_T)u_n(\mu_T))=\mu_T'(t)S(-\mu_T)(\partial_t+|\partial_x|\partial_x) u_n \big|_{\mu_T}.
$$ 
Thus, $\partial_t(S(-\mu_T)u_n(\mu_T))=0$ on $[-2T,0)$ since $\mu_T'(t)=0$ there. 
In addition, on $[0,T]$ we have $\partial_t(S(-\mu_T)u_n(\mu_T))=S(-t)(\partial_t+|\partial_x|\partial_x) u_n(t)$ while on $(T, 2T]$ we have $\partial_t(S(-\mu_T)u_n(\mu_T))=-S(t-2T)(\partial_t+|\partial_x|\partial_x) u_n(2T-t)$. 
Therefore,
$$
\no{u_n^*}_{X^{s-1,1}_{\tau=-|\xi|\xi}} 
\lesssim
\no{u_n}_{C([0,T];H^s(\mathbb T))}
+
\no{u_n}_{X^{s-1,1}_{\tau=-|\xi|\xi,T}}
+
\no{(\partial_t+|\partial_x|\partial_x) u_n}_{L^2([0,T];H^{s-1}(\mathbb T))}.
$$
To handle the third term, let $u_n^{**}\in X^{s-1,1}_{\tau=-|\xi|\xi}$ be any extension of $u_n \in X^{s-1,1}_{\tau=-|\xi|\xi, T}$. Then, 
\ddd{
\no{(\partial_t+|\partial_x|\partial_x) u_n}_{L^2([0,T];H^{s-1}(\mathbb T))}
&=
\no{(\partial_t+|\partial_x|\partial_x) u_n^{**}}_{L^2([0,T];H^{s-1}(\mathbb T))}
\nn\\
&=
\no{\mathcal{F}^{-1}((\tau+|\xi|\xi)\widetilde{u})}_{L^2([0,T];H^{s-1}(\mathbb T))}
\nn\\
&\lesssim
\no{\langle{\xi\rangle}^{s-1} \langle{\tau+|\xi|\xi\rangle} \widetilde{u_n^{**}}}_{L^2(\mathbb Z\times\mathbb R)}
=
\no{u_n^{**}}_{X^{s-1,1}_{\tau=-|\xi|\xi}}.
\nn
}
Hence, taking the infimum of this inequality over all extensions, we infer
$$
\no{u_n^*}_{X^{s-1,1}_{\tau=-|\xi|\xi}}
\lesssim
\no{u_n}_{C([0,T];H^s(\mathbb T))} 
+
\no{u_n}_{X^{s-1,1}_{\tau=-|\xi|\xi,T}}.
$$

In order to deduce inequality \eqref{ext-est2}  from the above inequality, it suffices to show that the left-hand side converges to $\no{u^*}_{X^{s-1,1}_{\tau=-|\xi|\xi}}$.
We have
$$
\no{u_n^*-u_m^*}_{X^{s-1,1}_{\tau=-|\xi|\xi}}
=
\no{(u_n-u_m)^*}_{X^{s-1,1}_{\tau=-|\xi|\xi}}
\lesssim
\no{u_n-u_m}_{C([0,T];H^s(\mathbb T))} 
+
\no{u_n-u_m}_{X^{s-1,1}_{\tau=-|\xi|\xi,T}}
$$
and, in addition, 
$\no{u_n^*-u_m^*}_{L^{\infty}(\mathbb R;H^s(\mathbb T))}
\lesssim
\no{u_n-u_m}_{C([0,T];H^s(\mathbb T))}$.
Therefore, $u_n^*$ is Cauchy in $X^{s-1,1}_{\tau=-|\xi|\xi}$ and $L^{\infty}(\mathbb R;H^s(\mathbb T))$ and has limits $v_1$ and $v_2$, respectively. Moreover, since 
\ddd{
&\no{u_n^*-v_1}_{L^{\infty}(\mathbb R;H^{s-1}(\mathbb T))}
\lesssim
\no{u_n^*-v_1}_{X^{s-1,1}_{\tau=-|\xi|\xi}},
\nn\\
&\no{u_n^*-v_2}_{L^{\infty}(\mathbb R;H^{s-1}(\mathbb T))}
\leqslant
\no{u_n^*-v_2}_{L^{\infty}(\mathbb R;H^s(\mathbb T))},
\nn
}
we infer that $u_n^*$ converges to both $v_1$ and $v_2$ in  $L^{\infty}(\mathbb R;H^{s-1}(\mathbb T))$ and hence $v_1=v_2$. Finally, since for any $u\in C([0,T];H^s(\mathbb T))$ we have
$\no{u^*}_{L^{\infty}(\mathbb R;H^s(\mathbb T))}
\lesssim
\no{u}_{C([0,T];H^s(\mathbb T))}$,
it follows that
$
\no{u_n^*-u^*}_{L^{\infty}(\mathbb R;H^s(\mathbb T))}
\lesssim
\no{u_n-u}_{C([0,T];H^s(\mathbb T))}
$
and, therefore, $u_n^*\rightarrow u^*$ in $C([0,T]; H^s(\mathbb T))$ and  in $X^{s-1,1}_{\tau=-|\xi|\xi, T}$, proving inequality \eqref{ext-est2}.
\end{Proof}

Back to \eqref{nls-temp1}, using the extensions $v^*$ and $u^*$  defined by \eqref{v-ext-def} and  \eqref{u-ext-def} we have
\ddd{
&\quad
\no{\eta_T \int_{t'=0}^t e^{i(t-t')\partial_x^2} \, \partial_x\Pi^+(\partial_x^{-1}v\cdot\Pi^-(u_x))~dt'}_{X^{s+a,\frac 12+\delta}_T}
\label{extensions-ineq}
\\
&\leqslant
\no{\eta_T \int_{t'=0}^t e^{i(t-t')\partial_x^2} \, \partial_x\Pi^+(\partial_x^{-1}v^*\cdot\Pi^-(u^*_x))~dt'}_{X^{s+a,\frac 12+\delta}}
\lesssim
\no{\partial_x\Pi^+(\partial_x^{-1}v^*\cdot\Pi^-(u^*_x))}_{X^{s+a,-\frac 12+\delta}}
\nn
}
with the second inequality due to the following well-known result (see, for example, Proposition 2.12 in \cite{t2006}):
\eee{\label{tao1-l}
\no{\eta(t)\int_{t'=0}^t e^{i(t-t')\partial_x^2} F(x,t') dt'}_{X^{s,b}}
\lesssim
\no{F}_{X^{s,b-1}}, \quad s\in\mathbb R, \ b>\tfrac 12.
}
We shall now estimate the right-hand side of \eqref{extensions-ineq}. Applying projections, we have
\ddd{
\no{\partial_x\Pi^+(\partial_x^{-1}v^*\cdot\Pi^-(u^*_x))}_{X^{s+a,-\frac 12+\delta}}
&\leqslant
\sum_{k=1}^\infty \sum_{j=1}^{k} \sum_{m=1}^{k}
\no{\partial_x P_j\Pi^+\big((\partial_x^{-1} v^*)_k\Pi^-(u_x^*)_m\big)}_{X^{s+a,-\frac 12+\delta}}
\label{lp-temp1}
\\
&\lesssim
\sum_{k=1}^\infty \sum_{j=1}^{k} \sum_{m=1}^{k}
2^{j(s+a+1)}
\no{P_j\Pi^+\big((\partial_x^{-1} v^*)_k\Pi^-(u_x^*)_m\big)}_{X^{0,-\frac 12+\delta}}
\nn
}
with the restrictions on the summation ranges due to the support properties of $\Pi^\pm$. 
Furthermore, employing the bilinear estimate \eqref{interp-c} we find
\ddd{
\no{P_j\Pi^+((\partial_x^{-1} v^*)_k\Pi^-(u_x^*)_m)}_{X^{0,-\frac 12+\delta}}
&\lesssim
\left[
2^{-\frac{m+j}{2} + k \left(\frac 16+\delta\right)} 
\left(\no{\Pi^-(u_x^*)_m}_{L^\infty(\mathbb R_t; L^2(\mathbb T))}
\right.\right.
\nn\\
&\quad \
\left.\left.
+
2^{-\frac{m+j}{2}} 
\no{\Pi^-(u_x^*)_m}_{\bar X^{0,1}}
\right)
\right]^{\theta}
\no{\Pi^-(u_x^*)_m}_{\bar X^{0,1}}^{1-\theta} \no{(\partial_x^{-1} v^*)_k}_{X^{0,\frac 12}}.
\nn
}
Moreover, noting that $\xi^2 = -|\xi|\xi$ for $\xi<0$, we have
$$
\no{\Pi^-(u_x^*)_m}_{\bar X^{0,1}}
\lesssim
2^{m(1-(\sigma-1))} \no{\Pi^-(u^*)_m}_{\bar X^{\sigma-1,1}}
\lesssim
2^{m(2-\sigma)} \no{\Pi^-(u^*)_m}_{X_{\tau = -|\xi|\xi}^{\sigma-1,1}}
$$
where $\sigma := \min\left\{s,\f{1}{2}\right\}$.  We denote $Z_T:= C([0,T];H^{\sigma}(\mathbb T))\cap X_{\tau=-|\xi|\xi, T}^{\sigma-1,1}$.  Then, employing Lemma \ref{u-ext-def-l}, we obtain
\ddd{
&\quad
\no{P_j\Pi^+((\partial_x^{-1} v^*)_k\Pi^-(u_x^*)_m)}_{X^{0,-\frac 12+\delta}}
\nn\\
&\lesssim
\left[
2^{-\frac{m+j}{2} + k\left(\frac 16+\delta\right)}
\left(2^{m(1-\sigma)}
+
2^{m(\frac{3}{2}-\sigma)-\frac{j}{2}}
\right)
\right]^{\theta}
2^{m(2-\sigma)(1-\theta)-k(1+s)}
\no{u}_{Z_T} \no{v}_{X^{s,\frac 12}_T}.
\nn
}
In turn, \eqref{lp-temp1} becomes
\ddd{
&\quad
\no{\partial_x\Pi^+(\partial_x^{-1}v^*\cdot\Pi^-(u^*_x))}_{X^{s+a,-\frac 12+\delta}}
\label{lp-temp2}
\\
&\lesssim
\no{u}_{Z_T} \no{v}_{X^{s,\frac 12}_T}
\sum_{k=1}^\infty \sum_{j=1}^{k} \sum_{m=1}^{k}
2^{j(s+a+1)+m(2-\sigma)(1-\theta)-k(1+s)}
\left[
2^{-\frac{m+j}{2} + k\left(\frac 16+\delta\right)}
\left(2^{m(1-\sigma)}
+
2^{m(\frac{3}{2}-\sigma)-\frac{j}{2}}
\right)
\right]^{\theta}
\nn
}
and, therefore, it suffices to control the multiplier
$$
M:=
2^{j(s+a+1)+m(2-\sigma)(1-\theta)-k(1+s)}
\left[
2^{-\frac{m+j}{2} + k\left(\frac 16+\delta\right)}
\left(2^{m(1-\sigma)}
+
2^{m(\frac{3}{2}-\sigma)-\frac{j}{2}}
\right)
\right]^{\theta}
$$
for $k, j, m$ as in \eqref{lp-temp2}.
Recalling that $\theta=\frac{1/2 - \delta}{1/2 + \delta}$ and $0<\delta \ll 1$, we may write $\theta=1-\epsilon$ for $0<\epsilon := \frac{2\delta}{1/2+\delta} \ll 1$. Then,  
\ddd{
M 
&=
2^{j(s+a+1)} \, 2^{-k(1+s)} \, 2^{m(2-\sigma)\epsilon} 
\left(2^{m(\frac 12-\sigma)} \, 2^{k(\frac 16+\delta)} \, 2^{-\frac{j}{2}}
+
2^{m(1-\sigma)} \, 2^{k(\frac 16+\delta)} \, 2^{-j}\right)^{1-\epsilon}
\nn\\
&=
\frac{
2^{j(s+a)} \, 2^{k(-\frac 56-s+\delta)} \, 2^{m(2-\sigma)\varepsilon} \, 2^{m(\frac 12-\sigma)}
\left(
2^{\frac j2} 
+
2^{\frac m2} \right)
}
{
\left(2^{m(\frac 12-\sigma)} \, 2^{k(\frac 16+\delta)} \, 2^{-\frac{j}{2}}
+
2^{m(1-\sigma)} \, 2^{k(\frac 16+\delta)} \, 2^{-j}\right)^\epsilon
}.
\nn
}
Hence, since $\epsilon>0$  and $2^{k(\frac 16+\delta)} > 1$, we have
$$
M \leqslant
2^{j(s+a+\frac \epsilon 2)} \, 2^{k(-\frac 56-s+\delta)} \, 2^{m (\frac 12-\sigma+\frac{3\epsilon}{2})}
\big(2^{\frac j2} + 2^{\frac m2} \big).
$$
Moreover, since $j,m\leqslant k$, if $s\leqslant \frac 12$ then $\sigma = s$, so we obtain
$$
M \lesssim
2^{k(s+a+\frac \epsilon 2)} \, 2^{k(-\frac 56-s+\delta)} \, 2^{k(\frac 12-s+\frac{3\epsilon}{2})} \, 2^{\frac k2}  
=
2^{k\left(a+2\epsilon+\frac 16+\delta-s\right)}
\leqslant
2^{k\left(a+\frac 16+9\delta-s\right)}
$$
while if $s>\frac 12$ then $\sigma = \f{1}{2}$, so we have
$$
M \leqslant
2^{j(s+a+\frac \epsilon 2)} \, 2^{k(-\frac 56-s+\delta)} \, 2^{m\frac{3\epsilon}{2}}
\big(2^{\frac j2} + 2^{\frac m2} \big)
\lesssim
2^{k \left(a+2\epsilon-\frac 13+\delta\right)}
\leqslant
2^{k \left(a-\frac 13+9\delta\right)}.
$$
Therefore, returning to \eqref{lp-temp2}, for $0\leqslant s\leqslant \frac 12$ we deduce
\eee{
\no{\partial_x\Pi^+(\partial_x^{-1}v^*\cdot \Pi^-(u^*_x))}_{X^{s+a,-\frac 12+\delta}}
\lesssim
\no{u}_{Z_T}  \no{v}_{X^{s,\frac 12}_T}
\sum_{k=1}^\infty \sum_{j=1}^{k} \sum_{m=1}^{k}
2^{k(a+\frac 16+9\delta-s)},
\nn
}
where the sum converges for  $a < s-\frac 16-9\delta$, while for $\frac 12\leqslant s\leqslant 1$ we deduce
\eee{
\no{\partial_x\Pi^+(\partial_x^{-1}v^*\cdot\Pi^-(u^*_x))}_{X^{s+a,-\frac 12+\delta}}
\lesssim
\no{u}_{Z_T}  \no{v}_{X^{s,\frac 12}_T}
\sum_{k=1}^\infty \sum_{j=1}^{k} \sum_{m=1}^{k}
2^{k(a-\frac{1}{3}+9\delta)},
\nn
}
where the sum converges for  $a<\frac{1}{3}-9\delta$.
The last two inequalities combined with inequality \eqref{extensions-ineq} and  definition \eqref{v-def} yield  the following bound for the right-hand side of \eqref{nls-temp1}:
\eee{\label{nl-smooth-bound-0}
\no{e^{-iKt} w - e^{it\partial_x^2}w_0}_{C([0,T];H^{s+a}(\mathbb T))}\lesssim 
\no{u}_{Z_T}  \no{w}_{X^{s,\frac 12}_T}
}
with $0<a<\min\left\{ s-\frac 16-9\delta, \frac{1}{3}-9\delta\right\}$, where we have used the fact that $\no{v}_{X^{s,\frac 12}_T} 
\lesssim
\no{w}_{X^{s,\frac 12}_T}$
since 
$\left\langle \tau - K +\xi^2 \right\rangle
\leqslant
\left(1 +  \left|\tau +\xi^2 \right|\right)\left(1+K\right)$.
Combining \eqref{nl-smooth-bound-0} with the estimates for $u$ and $w$ provided by Theorem \ref{bo-wp-t}, we conclude that
%
\eee{\label{nl-smooth-bound}
\no{e^{-iKt} w - e^{it\partial_x^2}w_0}_{C([0,T];H^{s+a}(\mathbb T))}
\lesssim 
\max\left\{\no{u_0}_{L^2(\mathbb T)}^{2s}, 1 \right\} \left(\no{u_0}_{H^{\sigma}(\mathbb T)}+\no{u_0}_{H^{\sigma}(\mathbb T)}^2 \right) \n{u_0}{H^s(\mathbb T)}
}
completing the proof of Theorem \ref{nl-smooth-t}.

\section{Polynomial Bound: Proof of Theorem \ref{poly-bound-t}}
\label{poly-bound-s}

We will now employ the nonlinear smoothing effect of Theorem \ref{nl-smooth-t}  in order to establish the polynomial bound of Theorem \ref{poly-bound-t}. 
We begin by noting that estimate \eqref{nl-smooth-bound} (which is the concrete expression of the nonlinear smoothing effect) for $s=\frac 12$ and $0< a < \frac 13$  implies 
\eee{\label{nl-smooth-pb}
\n{w(t) - e^{it(\p_x^2 + K)} w_0}{H^{\f{5}{6} - \ve}(\mathbb T)} \leqslant C\big(\n{u_0}{H^{\frac 12}(\mathbb T)}\big), \quad \ve := \tfrac 13-a>0, \ t\in [0, T],
}
where $C\big(\n{u_0}{H^{\frac 12}(\mathbb T)}\big)$ is a constant that depends only on $\n{u_0}{H^{\frac 12}(\mathbb T)}$. 

We also note that 
\eee{\label{wu-est}
\n{w(t)}{H^s(\mathbb T)} \leqslant C\big(s, \no{u_0}_{L^2(\mathbb T)}\big) \n{u(t)}{H^s(\mathbb T)}, \quad 0\leqslant s \leqslant 1, \ t\in \mathbb R.
}
Indeed, for $\frac 12 < s \leqslant 1$ inequality \eqref{wu-est} follows from the algebra property after recalling that $w \simeq \Pi^+(u e^{-iF/2})$ and observing that $\no{e^{-iF/2}}_{H^s(\mathbb T)} \leqslant \no{e^{-iF/2}}_{H^1(\mathbb T)} \leqslant  1+\no{u_0}_{L^2(\mathbb T)}$ from the physical definition of the $H^1$-norm and the conservation of the $L^2$-norm. 
Moreover, for $0\leqslant s \leqslant \frac 12$ inequality \eqref{wu-est} follows directly from inequality (2.13) of \cite{mp2012}. 

In addition,  the $H^{\frac 12}$-norm of $u$ can be controlled via  the following result.
\begin{lemma}\label{conservation-l}
Let $u$ satisfy the BO initial value problem \eqref{bo-ivp}. Then, 
\eee{\label{1/2-ineq}
\n{u(t)}{H^{\f{1}{2}}(\mathbb T)} \leqslant C\big(\n{u_0}{H^{\f{1}{2}}(\mathbb T)}\big), \quad t\in \mathbb R,
}
where $C\big(\n{u_0}{H^{\frac 12}(\mathbb T)}\big)$ is a constant that depends only on $\n{u_0}{H^{\frac 12}(\mathbb T)}$.
\end{lemma}

\begin{Proof}[Proof of Lemma \ref{conservation-l}]
Multiplying the BO equation \eqref{bo-eq} by $|\p_x|u$, which is defined via  Fourier transform by $\widehat{|\p_x|^j u}(\xi) = |\xi|^j \widehat u(\xi)$, and integrating  over $\mathbb T$, we have
\eee{\label{1/2-temp1}
\int_{x\in \mathbb T} u_t \cdot |\p_x| u \, dx + \int_{x\in \mathbb T} \mathcal{H} u_{xx} \cdot |\p_x| u \, dx = \frac 12 \int_{x\in \mathbb T} \p_x (u^2) \cdot  |\p_x| u\, dx.
}
For the  first integral, recalling that $u$ is real-valued, and hence that $\overline{\widehat u(\xi)} = \widehat u(-\xi)$ and in turn $|\p_x u| = \overline{|\p_x|u}$,  and using Parseval's identity twice, we find that
\eee{\label{1/2-temp2}
\int_{x\in \mathbb T} u_t \cdot |\p_x| u \, dx
=
\frac 12 \cdot \frac{1}{2\pi} \sum_{\xi\in\mathbb Z} \Big(\widehat u_t \overline{\widehat{|\p_x| u}} + \overline{\widehat u}_t \widehat{|\p_x| u} \Big)
=
\frac 12 \p_t \no{u}_{\dot H^\frac 12(\mathbb T)}^2,
}
where $\dot H^{\frac 12}$ denotes the homogeneous counterpart of $H^{\frac 12}$. 
Also, recalling in addition that $\mathcal H \p_x^2 = |\p_x| \p_x$ and using Parseval's identity, we find that the second integral vanishes:
\eee{\label{1/2-temp3}
\int_{x\in \mathbb T} \mathcal{H} u_{xx} \cdot |\p_x| u \, dx
=
\frac{1}{2\pi} \sum_{\xi \in \mathbb Z} \widehat{|\p_x|\p_x u} \cdot \overline{\widehat{|\p_x| u}}
=
\frac{1}{2\pi} \sum_{\xi \in \mathbb Z} \xi^3 \widehat u(\xi) \widehat u(-\xi) = 0.
}
Finally, integrating by parts and substituting from the BO equation, we write the third integral  as
\eee{\label{1/2-temp4}
\frac 12 \int_{x\in \mathbb T} \p_x (u^2) \cdot  |\p_x| u\, dx
=
-\frac 12 \int_{x\in \mathbb T} u^2 \cdot \mathcal H u_{xx}\, dx
=
\frac 16 \, \p_t \int_{x\in \mathbb T} u^3 \, dx.
}
Combining \eqref{1/2-temp1}-\eqref{1/2-temp4}, we deduce that the quantity $\no{u}_{\dot H^\frac 12(\mathbb T)}^2
-
\frac 13  \int_{x\in \mathbb T} u^3 \, dx$ is conserved, i.e. 
\eee{\label{1/2-temp5}
\no{u(t)}_{\dot H^\frac 12(\mathbb T)}^2
-
\frac 13 \int_{x\in \mathbb T} u^3(t) \, dx
=
\no{u_0}_{\dot H^\frac 12(\mathbb T)}^2
-
\frac 13 \int_{x\in \mathbb T} u_0^3 \, dx.
}
Moreover, by Sobolev's inequality (e.g. see Theorem 4.31 in \cite{af2003}), the fractional Sobolev-Gagliardo-Nirenberg inequality (see Corollary 1.5 in \cite{hmow2011})  and the   conservation of the $L^2$-norm, we have
\eee{\label{1/2-temp6}
\int_{x\in \mathbb T} u^3(t) \, dx
\leqslant
\no{u(t)}_{L^3(\mathbb T)}^3
\lesssim
\no{u(t)}_{H^{\frac 16}(\mathbb T)}^3
=
\no{u_0}_{L^2(\mathbb T)}^2 \no{u(t)}_{H^{\frac 12}(\mathbb T)}.
}
Combining \eqref{1/2-temp5} and \eqref{1/2-temp6}, we find
\ddd{
&\quad
\frac 13 \Big(\no{u(t)}_{\dot H^\frac 12(\mathbb T)}^2
-
 \no{u_0}_{L^2(\mathbb T)}^2 \no{u(t)}_{H^{\frac 12}(\mathbb T)} \Big)
\nn\\
&\lesssim
\no{u_0}_{\dot H^\frac 12(\mathbb T)}^2
-
\frac 13 \int_{x\in \mathbb T} u_0^3 \, dx
\leqslant
\no{u_0}_{\dot H^\frac 12(\mathbb T)}^2
+
\no{u_0}_{L^3(\mathbb T)}^3
\lesssim
\no{u_0}_{\dot H^\frac 12(\mathbb T)}^2
+
\no{u_0}_{L^2(\mathbb T)}^2 \no{u_0}_{H^{\frac 12}(\mathbb T)},
\nn
}
i.e.
\eee{\label{1/2-temp7}
\no{u(t)}_{\dot H^\frac 12(\mathbb T)}^2
-
 \no{u_0}_{L^2(\mathbb T)}^2 \no{u(t)}_{H^{\frac 12}(\mathbb T)}
\leqslant
C\big(\no{u_0}_{H^{\frac 12}(\mathbb T)}\big).
}
But note that for $\xi \in \mathbb Z \setminus \{0\}$ we have $|\xi| \simeq \left\langle \xi \right\rangle$. 
Using this fact together with our assumption of mean-zero data, we infer from \eqref{1/2-temp7} the inequality
\eee{
\no{u(t)}_{H^\frac 12(\mathbb T)}^2
-
 \no{u_0}_{L^2(\mathbb T)}^2 \no{u(t)}_{H^{\frac 12}(\mathbb T)}
\leqslant
C\big(\no{u_0}_{H^{\frac 12}(\mathbb T)}\big).
\nn
}
Completing the square on the left-hand side yields the desired inequality \eqref{1/2-ineq}.
\end{Proof}

Before proceeding to the proof of Theorem \ref{poly-bound-t}, we establish the following inequality.
\begin{proposition}\label{uw-p}
Suppose that $\frac 12<s\leqslant 1$. Then, for all $t\in \mathbb R$ we have
\eee{\label{uw-ineq}
\no{u(t)}_{H^s(\mathbb T)}
\lesssim
\big(1+\no{u_0}_{L^2(\mathbb T)}\big)
\left(\no{w(t)}_{H^s(\mathbb T)}+
\big(1+\no{u_0}_{L^2(\mathbb T)}\big)
\big[1+C\big(\no{u_0}_{H^{\frac 12}(\mathbb T)}\big)\big]\right).
}
\end{proposition}

\begin{Proof}[Proof of Proposition \ref{uw-p}]
We suppress the $t$-dependence for brevity. Note that $u=u^+ + \overline{u^+}$ so
$
\no{u}_{H^s(\mathbb T)}
\leqslant
2\no{u^+}_{H^s(\mathbb T)}
$.
Also, 
$
u=2i e^{iF/2} w + 2i e^{iF/2}\partial_x\Pi^-(e^{-iF/2})
$
and hence
$$
\no{u^+}_{H^s(\mathbb T)}
\lesssim
\big\|\Pi^+(e^{iF/2} w)\big\|_{H^s(\mathbb T)}
+
\big\|\Pi^+\big(e^{iF/2}\partial_x\Pi^-(e^{-iF/2})\big)\big\|_{H^s(\mathbb T)}.
$$
By Lemmas 3.1 and 3.2 of \cite{m2007}  we have  
$$
\big\|\Pi^+(e^{iF/2} w)\big\|_{H^s(\mathbb T)}
\lesssim
\no{w}_{H^s(\mathbb T)} \big(1+\no{u_0}_{L^2(\mathbb T)}\big).
$$
and, for $s_1+s_2=s+1$, $s_1\geqslant s$ and $s_2\geqslant 0$,
$$
\big\|\Pi^+(e^{iF/2}\partial_x\Pi^-(e^{-iF/2}))\big\|_{H^s(\mathbb T)}
\lesssim
\big\|J_x^{s_1} e^{iF/2}\big\|_{L^4(\mathbb T)}
\big\|J_x^{s_2} e^{-iF/2}\big\|_{L^4(\mathbb T)}.
$$
Since
$$
\big\|J_x^{s_1} e^{iF/2}\big\|_{L^4(\mathbb T)}
\lesssim
\big\|J_x^{s_1+\frac{1}{4}} e^{iF}\big\|_{L^2(\mathbb T)},
\quad
\big\|J_x^{s_2} e^{-iF/2}\big\|_{L^4(\mathbb T)}
\lesssim
\big\|J_x^{s_2+\frac{1}{4}} e^{-iF/2}\big\|_{L^2(\mathbb T)}
$$
by the Sobolev embedding, taking $s_2=\frac{3}{4}$  we have
$$
\big\|J_x^{s_2+\frac{1}{4}} e^{-iF/2}\big\|_{L^2(\mathbb T)}
\lesssim
\big\|e^{-iF/2}\big\|_{H^1(\mathbb T)}
\lesssim
\big(1+\no{u_0}_{L^2(\mathbb T)}\big).
$$
Then, $s_1=s+\frac{1}{4}$  and for $s=\frac 12+\delta$, $0 < \delta\leqslant\frac 12$, we find
\ddd{
\big\|J_x^{s_1+\frac{1}{4}} e^{iF/2}\big\|_{L^2(\mathbb T)}
=
\big\|J_x^{1+\delta} e^{iF/2}\big\|_{L^2(\mathbb T)}
&\lesssim
\big\|J_x^{\delta} e^{iF/2}\big\|_{L^2(\mathbb T)}
+
\big\|J_x^{\delta}(u e^{iF/2})\big\|_{L^2(\mathbb T)}
\nn\\
&\lesssim
(1+\no{u_0}_{L^2(\mathbb T)})+\big\|J_x^{\delta}u\big\|_{L^2(\mathbb T)} \big(1+\no{u_0}_{L^2(\mathbb T)}\big)
\nn
}
with  the second inequality due to Lemma 3.1 of \cite{m2007}. Noting further that   Lemma \ref{conservation-l} implies $\big\|J_x^{\delta}u\big\|_{L^2(\mathbb T)}
\lesssim
C(\no{u_0}_{H^{\frac 12}(\mathbb T)})$, we obtain the desired inequality \eqref{uw-ineq}.
\end{Proof}

We now combine inequalities \eqref{nl-smooth-pb}-\eqref{1/2-ineq} with inequality \eqref{uw-ineq} to obtain the polynomial bound of Theorem \ref{poly-bound-t}.
First, consider $\f{1}{2} < s < \f{5}{6}$.  Given $u_0 \in H^s$, let $T= T(\n{u_0}{L^2})$ be as in Theorem~\ref{bo-wp-t}.  Suppose $t \in [nT, (n+1)T)$ for some $n\in \mathbb N \cup \{0\}$.  Then, write
\eee{\label{bd-temp0}
w(t) = Q_{\leqslant n^3} w(t) + Q_{>n^3} w(t),
}
where $Q_{\leqslant n^3}$ and $Q_{>n^3}$ are the projections onto Fourier modes\footnote{Here, we use different notations for the projections, as $P_k=Q_{2^k}$.} 
  whose absolute value is less than or equal to $n^3$ and greater than $n^3$, respectively.  For the first component, we have
\eee{\label{bd-temp1}
\n{Q_{\leqslant n^3} w(t)}{H^s(\mathbb T)} \lesssim n^{3(s-\f{1}{2})} \n{w(t)}{H^{\f{1}{2}}(\mathbb T)} \lesssim \lan{n}^{3(s-\f{1}{2})} \n{u(t)}{H^{ \f{1}{2}}(\mathbb T)} \lesssim  \lan{t}^{3(s-\f{1}{2})}C\big(\n{u_0}{H^{\f{1}{2}}(\mathbb T)}\big).
}
where the final implicit constant depends on $T = T(\n{u_0}{L^2})$.
Hence, it remains to control the second component, which we rewrite as
$$
Q_{>n^3} \left(w(t) - e^{i(t-nT)(\p_x^2 + K)} w(nT)\right) + Q_{>n^3} e^{i(t-nT)(\p_x^2 + K)} w(nT).
$$
Since $s<\f{5}{6}$, employing estimate \eqref{nl-smooth-pb} after shifting the time interval from $[0, T]$ to $[nT, (n+1)T]$ together with estimate  \eqref{1/2-ineq}, we can control the first  part above as follows: 
\begin{align*}
&\ \quad
\n{Q_{>n^3} \left(w(t) - e^{i(t-nT)(\p_x^2 + K)} w(nT)\right)}{H^s(\mathbb T)}
\nn\\
&= \Big\|J_x^{s-\f{5}{6}+\ve} Q_{>n^3} J_x^{\f{5}{6}-\ve} \left(w(t) - e^{i(t-nT)(\p_x^2 + K)} w(nT)\right)\Big\|_{L^2(\mathbb T)}
\\
&\lesssim n^{3s-\f{5}{2}+3\ve} \n{ w(t) - e^{i(t-nT)(\p_x^2 + K)} w(nT) }{H^{\f{5}{6}-\ve}(\mathbb T)}
\\
&\lesssim n^{3s-\f{5}{2}+3\ve} C\big(\n{u(Tn)}{H^{\f{1}{2}}(\mathbb T)}\big) \lesssim n^{3s-\f{5}{2}+3\ve} C\big(\n{u_0}{H^{\f{1}{2}}(\mathbb T)}\big).
\end{align*}
For the second part,  writing 
$$
Q_{>(n-1)^3} w(nT) = Q_{>(n-1)^3}\Big(w(nT) - e^{iT (\p_x^2 + K)}  w((n-1)T)  \Big) + Q_{>(n-1)^3}  e^{iT (\p_x^2 + K)}  w((n-1)T)
$$
and the \textit{strict} inequality
$$
\n{Q_{>n^3} e^{i(t-nT)(\p_x^2 + K)} w(nT)}{H^s(\mathbb T)} \leqslant \n{ Q_{>(n-1)^3} w(nT)}{H^s(\mathbb T)}
$$
allow us to repeat our earlier computation for the first part to obtain
$$
\n{ Q_{>n^3} e^{i(t-nT)(\p_x^2 + K)} w(nT)}{H^s(\mathbb T)} \leqslant  \widetilde C \left(n-1\right)^{3s-\f{5}{2}+3\ve} C\big(\n{u_0}{H^{\f{1}{2}}(\mathbb T)}\big) +  \n{ Q_{>(n-2)^3} w((n-1)T)}{H^s(\mathbb T)}.
$$
As before, it is important that the second term on the right-hand side does not pick up any constant. Thus, we can iterate this process $n$ times to obtain
\eee{\label{bd-temp2}
\n{Q_{>n^3} w(t)}{H^s(\mathbb T)} \leqslant \sum_{k=1}^n k^{3s-\f{5}{2}+3\ve}  C\big(\n{u_0}{H^{\frac 12}(\mathbb T)}\big) + \n{w_0}{H^s(\mathbb T)} \lesssim n^{3(s - \f{1}{2}+\ve)}   C\big(\n{u_0}{H^{\frac 12}(\mathbb T)}\big) + \n{u_0}{H^s(\mathbb T)},
}
where the implicit constant in the second inequality depends on $s$ and $\ve$ and where we have used the following lemma.

\begin{lemma}\label{riemann-l}
For $\al>-1$ and $N\gg1$,
$$
\sum_{k=1}^N k^{\al} =\frac{1}{\alpha+1} N^{\alpha+1}+O\big(N^{\max(0, \alpha)}\big).
$$
In particular, $\sum_{k=1}^N k^{\al}\leqslant C_\al N^{\alpha+1}$. 
\end{lemma}

\begin{Proof}[Proof of Lemma \ref{riemann-l}]
We   have 
$$
\sum_{k=1}^N k^{\al} -\frac{1}{\alpha+1} N^{\alpha+1}=\sum_{k=1}^N \int_{x=k-1}^k k^\al dx - \int_{x=0}^N x^\al dx
= 
\sum_{k=1}^N \int_{x=k-1}^k \left(k^\al-x^\al\right) dx.
$$
Thus, noting that $\left|k^\al-x^\al\right|\lesssim k^{\al-1}$ for $x\in \left[k-1,k\right]$, for $-1<\alpha<0$ we get the bound 
$$
\left|\sum_{k=1}^N k^{\al} -\frac{1}{\alpha+1} N^{\alpha+1}\right| \lesssim  \sum_{k=1}^N k^{\al-1}\lesssim 1
$$
while for $\al\geqslant 0$ we have
$$
\left|\sum_{k=1}^N k^{\al} -\frac{1}{\alpha+1} N^{\alpha+1}\right| 
=
\sum_{k=1}^N \int_{x=k-1}^k \left(k^\al-x^\al \right) dx
\leqslant
\sum_{k=1}^N \int_{x=k-1}^k \left[k^\al-(k-1)^\al \right] dx
= 
N^\al.
$$
The proof of the lemma is complete.
\end{Proof}

Note that $3s - \f{5}{2} >-1 \Leftrightarrow s > \f{1}{2}$ and so Lemma \ref{riemann-l} can be employed to yield the last inequality in \eqref{bd-temp2}.
Overall, combining \eqref{bd-temp1} and \eqref{bd-temp2} with the decomposition \eqref{bd-temp0}, for any $n\in \mathbb N \cup \{0\}$ and $t\in [nT, (n+1)T)$ we obtain
\ddd{
\n{w(t)}{H^s(\mathbb T)} 
&\leqslant  \n{Q_{\leqslant n^3} w(t)}{H^s(\mathbb T)} + \n{Q_{>n^3} w(t)}{H^s(\mathbb T)} 
\nn\\
&\lesssim \lan{t}^{3(s-\f{1}{2}+\ve)}C(\n{u_0}{H^{\frac 12}(\mathbb T)})  + \n{u_0}{H^s(\mathbb T)} 
\lesssim \lan{t}^{3(s-\f{1}{2}+\ve)}C(\n{u_0}{H^{s}(\mathbb T)}), 
\nn
}
where the implicit constants depend on $s$, $T$ and $\ve$. Therefore, using inequality \eqref{uw-ineq} we obtain the desired bound, concluding the proof of Theorem \ref{poly-bound-t} for $\frac 12 < s < \frac 56$.

For $\f{5}{6} \leqslant s \leqslant 1$, we can follow a similar computation to establish the same polynomial-in-time bound. Indeed, as before, we write $w(t) = Q_{\leqslant n^3} w(t) + Q_{>n^3} w(t)$ and note that the first component can be estimated once again as in \eqref{bd-temp1}.
Furthermore, we rewrite the second component as
$$
Q_{>n^3} \left(w(t) - e^{i(t-nT)(\p_x^2 + K)} w(nT)\right) + Q_{>n^3} e^{i(t-nT)(\p_x^2 + K)} w(nT)
$$
and note that estimate \eqref{nl-smooth-bound}  with $s=\frac 56 - \ve$ and $a = \frac 13 -\ve$  after shifting $[0, T]$ to $[nT, (n+1)T]$ implies
\begin{align*}
&\ \quad\n{Q_{>n^3} \left(w(t) - e^{i(t-nT)(\p_x^2 + K)} w(nT)\right)}{H^s(\mathbb T)}
\nn\\
&= \Big\|J_x^{s-\f{7}{6}+2\ve} Q_{>n^3} J_x^{\f{7}{6}-2\ve} \left(w(t) - e^{i(t-nT)(\p_x^2 +  K)} w(nT)\right)\Big\|_{L^2(\mathbb T)}
\\
&\lesssim n^{3s-\f{7}{2}+6\ve} \n{ w(t) - e^{i(t-nT)(\p_x^2 + K)} w(nT) }{H^{\f{7}{6}-2\ve}(\mathbb T)}
\\
&\lesssim n^{3s-\f{7}{2}+6\ve} C\big(\n{u(Tn)}{H^{\f{1}{2}}(\mathbb T)}\big) \n{u(nT)}{H^{\f{5}{6} - \ve}(\mathbb T)}\\
&\lesssim n^{3s-\f{7}{2}+6\ve} C\big(\n{u_0}{H^{\f{1}{2}}(\mathbb T)}\big) \lan{nT}^{1- 2\ve}
 \lesssim n^{3s - \f{5}{2}+ 4\ve} C\big(T, \n{u_0}{H^{\f{1}{2}}(\mathbb T)}\big).
\end{align*}
where we have also employed the previously established polynomial bound for $\n{u(nT)}{H^{\f{5}{6} - \ve}(\mathbb T)}$ to obtain the penultimate inequality.
Hence, repeating the iterative procedure used in the case $\f{1}{2} \leqslant s \leqslant \f{5}{6}$, we obtain the desired bound.  This concludes proof of Theorem~\ref{poly-bound-t}.

\vspace*{5mm}

\noindent
Bradley Isom \hfill Dionyssios Mantzavinos
\\
Department of Mathematics \hfill Department of Mathematics \\
University of Kansas  \hfill University of Kansas
\\
Lawrence, KS 66045 \hfill Lawrence, KS 66045 \\
\textit{bsisom@ku.edu} \hfill \textit{mantzavinos@ku.edu}

\vspace*{5mm}

\noindent
Seungly Oh \hfill Atanas Stefanov
\\
Department of Mathematics \hfill Department of Mathematics \\
Western New England University  \hfill University of Kansas
\\
Springfield, MA 01119 \hfill Lawrence, KS 66045 \\
\textit{seungly.oh@wne.edu} \hfill  \textit{stefanov@ku.edu}

\end{document}